\documentclass[11pt]{article}
\usepackage{epsfig}
\usepackage{amssymb}
\newtheorem{theorem}{Theorem}[section]
\newtheorem{lemma}[theorem]{Lemma}
\newtheorem{corollary}[theorem]{Corollary}

\newtheorem{proposition}[theorem]{Proposition}
\newtheorem{definition}{Definition}[section]

\newcommand{\proof}{\medskip \noindent\bf Proof. \rm}

\newcommand{\qed}{\hbox{\hspace{5mm}} \hfill \rule{7pt}{7pt} \\ \medskip}

\def\ba{\begin{array}}
\def\ea{\end{array}}
\def\beq{\begin{equation}}
\def\eeq{\end{equation}}
\def\bea{\begin{eqnarray}}
\def\eea{\end{eqnarray}}
\def\beann{\begin{eqnarray*}}
\def\eeann{\end{eqnarray*}}

\def\square{\vbox{\hrule\hbox{\vrule\kern3pt
    \vbox{\kern6pt}\kern3pt\vrule}\hrule}}

\def\diag{\mbox{\rm diag}}

\title{\bf WHEN DOES THE POSITIVE SEMIDEFINITENESS
CONSTRAINT HELP IN LIFTING PROCEDURES\thanks{Parts of the final version of this
paper were written while both authors were members of the Fields
Institute, Toronto, during Fall 1999.}}
\author{
Michel X. Goemans\thanks{MIT, Dept.~of Mathematics, Room 2-351,
Cambridge, MA 02139. Part of this research was performed when this
author was visiting the Department of Combinatorics and Optimization
of the University of Waterloo, and their hospitality and support are
gratefully acknowledged.  Research of this author was also supported
in part by NSF contract 9623859-CCR. }
\and
Levent Tun{\c c}el\thanks{
Department of Combinatorics and Optimization,
Faculty of Mathematics, University of Waterloo,
Waterloo, Ontario N2L 3G1, Canada.
Research of this author was supported
in part by a grant from NSERC and a PREA from
Ontario, Canada.}
}

\begin{document}
\maketitle
\begin{abstract}
We study the lift-and-project procedures of Lov\'asz and Schrijver for
0-1 integer programming problems. We prove that the procedure using
the positive semidefiniteness constraint is not better than the
one without it, in the worst case.
Various examples are considered.  We
also provide geometric conditions characterizing when the positive
semidefiniteness constraint does not help. 
\end{abstract}

\pagestyle{myheadings}
\thispagestyle{plain}
\markboth{GOEMANS and TUN{\c C}EL}
{EFFICIENCY OF LIFTING PROCEDURES}

\noindent
{\bf Keywords:}
Semidefinite lifting, semidefinite
programming, lift-and-project, integer programming

\noindent
{\bf AMS Subject Classification:} 90C10, 90C27, 47D20

\newpage

\section{Introduction}

Lov{\'a}sz and Schrijver (1991) have proposed a very intriguing
successive convex relaxation procedure for 0-1 integer programming
problems.  The procedure called $N_+$, to be defined shortly, when
applied to a classical linear programming (LP) relaxation of the
stable set problem (with only the edge and nonnegativity constraints)
produces a relaxation for which many well-known inequalities are
valid, including the odd hole, odd antihole, odd wheel, clique, and
even the orthonormal representation inequalities of Gr\"otschel,
Lov\'asz and Schrijver (1981). This implies that for many classes of
graphs, including perfect (for which clique inequalities are
sufficient) or t-perfect graphs (for which odd hole inequalities are
sufficient), one can find the maximum stable set by using the $N_+$
procedure.

The $N_+$ procedure is a strengthening of another procedure, called
$N$, also introduced by Lov\'asz and Schrijver. The main difference between
the two procedures is that $N_+$ involves a positive semidefinite
constraint. When applied to a linear programming relaxation, $N$ will
produce another (stronger) LP relaxation while $N_+$ will produce a
semidefinite relaxation. For the stable set problem,
Lov\'asz and Schrijver have shown that the relaxation produced by $N$
is much weaker than the one derived from $N_+$. 

In general, it is however not clear in which situations the procedure
$N_+$ is better or significantly better than $N$; especially, when
$N$ and $N_+$ are applied iteratively. In this paper, we
try to shed some light on this question. We generalize certain
properties derived by Lov\'asz and Schrijver. We also identify certain
situations in which $N$ produces the same relaxation as $N_+$. Several
examples are discussed throughout the paper, including one in which
the number of iterations of the $N_+$ procedure needed to derive the
convex hull of 0-1 points is equal to the dimension of the space,
hence resolving a question left open by Lov\'asz and Schrijver.

In the next section, we review the lift-and-project procedures
and their basic properties. Section \ref{sec:up} includes upper bounds
on the number of major iterations required by such procedures.
Section \ref{sec:low} discusses techniques to prove lower bounds
on the number of major iterations required. Sections
\ref{sec:lowercomp} and \ref{sec:gen} include 
geometric properties and characterizations of the
convex relaxations produced by the procedures. 

\section{Lov{\'a}sz-Schrijver procedures $N$ and $N_+$} \label{sec:n}

First, we describe two lift-and-project procedures proposed by Lov{\'
a}sz and Schrijver (1991) which produce tighter and tighter
relaxations of the convex hull of $0$-$1$ points in a convex set.  In
what follows, $e_j$ is the $j$th unit vector and $e$ is the vector of
all ones. The sizes of $e$ and $e_j$ will be clear from the context.
The cone generated by all 0-1 vectors $x \in {\mathbb R}^{d+1}$ with $x_0
= 1$ is called $Q$.  Let $K \subset Q$ denote a convex cone; for
example, $K$ could be a polyhedral cone obtained from a polytope $P$
in $[0,1]^d$ via homogenization using a new variable $x_0$. That
is, if
\[
P = \{x \in {\mathbb R}^d: \,\, Ax \leq b, \,\, 0 \leq x \leq e \},
\]
then
\[
K := \left\{ \pmatrix{x_0 \\ \cr
x \cr} \in {\mathbb R}^{d+1}:
\,\, Ax \leq x_0 b, \,\, 0 \leq x \leq x_0 e \,\,
\right\}.
\]
We are interested in determining (or approximating) $K_I$, the cone
generated by all 0-1 vectors of $K$.

Let $K^*$, $Q^*$ denote the dual cones
of $K$ and $Q$ under the standard Euclidean 
inner-product, e.g.,
$$K^* := \{s \in {\mathbb R}^{d+1}: x^Ts \,\, \geq \,\, 0, \,\, \forall
\, x \in K\}.$$ ${\cal S}^{d+1}$ denotes the space of $(d + 1)
\times(d+1)$ symmetric matrices and ${\cal S}_+^{d+1}$ denotes the
cone of $(d + 1) \times(d+1)$ symmetric, positive semidefinite
matrices. For a matrix $A\in {\cal S}^{d+1}$, we denote its positive
semidefiniteness by $A\succeq 0$.  When we deal with the
duals of convex cones in the space of $(d+1) \times (d+1)$ matrices
(or in the subspace of the symmetric matrices), we always
take the underlying inner-product to be the trace inner-product (or
Frobenius inner-product): $\langle A,B\rangle := Tr(A^TB)$.

Let $\mbox{diag}: {\cal S}^{d+1} \to {\mathbb R}^{d+1}$
denote the linear operator which maps a symmetric matrix to its
diagonal.
Then its adjoint $\mbox{diag}^* : {\mathbb R}^{d+1}
\to {\cal S}^{d+1}$ is the linear operator $\mbox{Diag}(\cdot)$
which maps a vector from ${\mathbb R}^{d+1}$ to the diagonal

atrix in ${\cal S}^{d+1}$ whose $(i,i)$th component is
the $i$th component of the original vector.

\begin{definition}[Lov\'asz and Schrijver (1991)]
\label{defi:M}
A $(d+1) \times (d+1)$ symmetric matrix, $Y$, with real entries is in
$M(K)$ if
\begin{itemize}
\item[(i)] $Ye_0 = \mbox{diag}(Y)$, and
\item[(ii)] $u^T Y v \, \geq \, 0, \,\,\, \forall \, u \in Q^*, \, v \in K^*.$
\end{itemize}
\end{definition}
Lov{\' a}sz and Schrijver note that condition (ii) of the above definition
is equivalent to $Y Q^* \subseteq K$ (where $YQ^*=\{Yx: x\in Q^*\}$),
or: $(ii)^{'}$ $Ye_i \in K$ for all $i \in \{1, \ldots, d\}$ and
$Y(e_0-e_i) \in K$ for all $i \in \{1, \ldots, d\},$ since the extreme
rays (after normalization) of the cone $Q^*$ are given by
$\mbox{ext}(Q^*) = \{e_1, e_2, \ldots, e_d, (e_0-e_1), (e_0-e_2),
\ldots, (e_0-e_d)\}.$

\begin{definition}[Lov\'asz and Schrijver (1991)]
\label{defi:Mplus}
$Y \in M_+(K)$ if $Y \in M(K)$ and $Y$ is positive semidefinite.
\end{definition}
Observe that if we take any $x\in K$ (not necessarily integral) and
consider $Y=xx^T$, $Y$ satisfies $Y\succeq 0$ and also (ii)', but this
specific $Y$
satisfies (i) if and only if $x$ is such that $x_i(x_0-x_i)=0$ for all
$i$, i.e.\ $x$ corresponds to a $0$-$1$ vector.

Now, we define the projections of these liftings $M$ and $M_+$:
$$N(K) := \{ \mbox{diag}(Y): Y \in M(K)\},$$
$$N_+(K) := \{ \mbox{diag}(Y): Y \in M_+(K)\}.$$ The above argument
regarding $xx^T$ shows that $K_I\subseteq N_+(K) \subseteq N(K)
\subseteq K$, the last inclusion following from the fact that
$Y(e_0-e_i)\in K$ and  $Ye_i \in K$ imply that $x=Ye_0\in K$. 

If $P$ is a polytope (or any convex set) in $[0,1]^d$ then we simply
write $N_+(P)$ to represent $\left\{x: \pmatrix{1 \cr x} \in
N_+(K)\right\}$ where $K$ is the cone obtained via homogenization
using the variable $x_0$, and similarly for $N(P)$. We also let
$M(P)=M(K)$ and $M_+(P)=M_+(K)$. 

\begin{figure}[htbp]
\hspace*{\fill}  \epsfig{figure=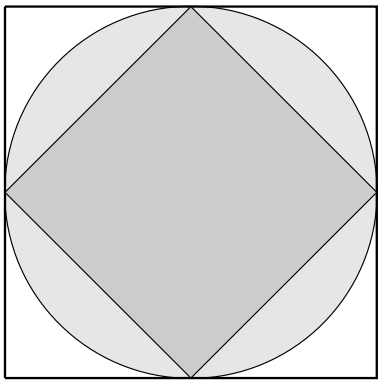,width=2in} \hspace*{\fill}

\caption{\label{fig1} Two convex sets with the same $N_+(\cdot)$.}
\end{figure}

We should point out that the definition of $M$ (or $M_+$) is such
that $M(K)$ depends only on the sets $K\cap \{x: x_i=x_0\}$ and $K\cap
\{x: x_i=0\}$ for all $i$. In particular, we have:
\begin{lemma}
Let $K$ and $K'$ be such that $K\cap \{x: x_i=x_0\}=K'\cap \{x:
x_i=x_0\}$ and $K\cap \{x: x_i=0\}= K'\cap \{x: x_i=0\}$ for all
$i \in \{1,\ldots,d\}$. Then $M(K)=M(K')$ (and $N(K)=N(K')$) and
$M_+(K)=M_+(K')$ (and $N_+(K)=N_+(K')$). 
\end{lemma}
For example,  $P=\{x\in {\mathbb R}^2: ||x-0.5e||_2 \leq
\frac{1}{2} \}$ and $P'=\{x\in {\mathbb R}^2: ||x-0.5e||_1 \leq
0.5\}$ (see Figure \ref{fig1}) have the same $N(P)=N(P')$.

The definitions of $M$, $N$, $M_+$ and $N_+$ are invariant under
various operations including flipping coordinates $x_i
\rightarrow (1-x_i)$ for any subset of the indices $\{1,2, \ldots,
d\}$. More formally,
\begin{proposition}[Lov\'asz and Schrijver (1991)]
\label{prop:2.1}
Let $A$ be a linear transformation mapping $Q$ onto itself. 
Then
\[
N(AK) = AN(K)
\mbox{ and } 
N_+(AK) = AN_+(K).
\]
\end{proposition}

One crucial feature of the operators $N$ and $N_+$ is that they can be
iterated.  The iterated operators $N^r(K)$ and
$N_+^r(K)$ are defined as follows.
$N^0 (K) := K$, $N_+^0 (K) := K$, $N^r(K) := N (N^{r-1}(K))$
and $N_+^r(K) := N_+ (N_+^{r-1}(K))$
for all integers $r \geq 1$. Lov\'asz and Schrijver (1991) show
that, even without the positive semidefiniteness
constraints, $d$ iterations are
sufficient to get $K_I$:
\begin{theorem}[Lov\'asz and Schrijver (1991)]
\label{thm:LS}
$$ K \supseteq N(K) \supseteq N^2(K) \supseteq \ldots
\supseteq N^d(K) = K_I$$
and
$$ K \supseteq N_+(K) \supseteq N_+^2(K) \supseteq \ldots
\supseteq N_+^d(K) = K_I.$$
\end{theorem}

Let $a^T x \leq \alpha x_0$ be a valid inequality for $K_I$. Then the
smallest nonnegative integer $r$ such that $a^T x \leq \alpha x_0$ is
valid for $N^r(K)$ is called {\em the $N$-rank of $a^T x \leq \alpha
x_0$ relative to $K$.}  The $N_+$-rank of $a^T x \leq \alpha x_0$
relative to $K$ is defined similarly. The above theorem states that
these ranks are at most $d$ for any valid inequality. The $N$-rank (resp.\
$N_+$-rank) of a cone $K$
is the smallest nonnegative integer $r$ such that $N^r(K)=K_I$ (resp.\
$N^r_+(K)=K_I$).

Theorem \ref{thm:LS} can also be proved using the results of Balas
(1974), see Balas, Ceria and Cornu\'ejols (1993).  Our interest, in
this paper, mostly lies in understanding the strength of $N_+$ in
comparison to $N$. Consider the stable set polytope on a graph
$G=(V,E)$ defined as the convex hull of incidence vectors of sets of
non-adjacent vertices (known as stable sets).  Let $FRAC$ be the relaxation
defined by the edge constraints ($x_i+x_j\leq 1$ for all edges
$(i,j)\in E$) and the nonnegativity constraints ($x_i\geq 0$ for all
$i\in V$).  Then $N(FRAC)$ is exactly equal to the relaxation obtained by
adding all odd hole inequalities, saying that $\sum_{i\in C} x_i \leq
\frac{|C|-1}{2}$ for any odd cycle $C$ with no chords. However,
many more complicated inequalities have small $N_+$-rank.  Lov\'asz
and Schrijver (1991) prove that odd hole, odd antihole, odd wheel,
clique and orthogonal inequalities all have $N_+$-rank at most 1,
relative to $FRAC$. These results are proved using Lemma \ref{lem:2.2}
of next section, except for the orthogonality constraints. In
contrast, the $N$-rank of a clique inequality for example is equal to
$p-2$ where $p$ is the size of the clique.  Note that the separation
problem for the class of clique inequalities is NP-hard (and so is the
problem of optimizing over the clique inequalities, see Gr\"otschel,
Lov\'asz and Schrijver (1981)). $N_+$, however, leads to a
polynomial-time separation algorithm for a broader class of
inequalities. This, and more generally the importance of $N$ and
$N_+$, stems from the following result. 

\begin{theorem}[Lov\'asz and Schrijver (1991)]
\label{thm:complexity}
If we have a weak separation oracle for $K$ then we have a weak
separation oracle for $N^r(K)$ and $N_+^r(K)$ for any fixed constant $r$.
\end{theorem}
Together with the equivalence between (weak) optimization and (weak)
separation (Gr\"otschel et al.~(1981)), this implies for example that
the stable set problem can be solved in polynomial time for any graph
with bounded $N_+$-rank (Lov\'asz and Schrijver~(1991)). 

Next we study the upper bounds on $N$- and $N_+$-ranks of
inequalities and convex sets.

\section{Upper bounds on the $N$- and $N_+$-rank} \label{sec:up}

Lov\'asz and Schrijver give some ways to upper bound the $N$-rank of
an inequality. They show the following. 
\begin{lemma}[Lov\'asz and Schrijver (1991)]
\label{lem:2.1}
\[
N_+(K) \subseteq N(K)  \subseteq \left(K \cap \{x:  x_i = 0
\}\right) + \left(K \cap \{x: x_i = x_0 \}
\right),
\mbox{ for all } i \in \{ 1,2, \ldots, d \}.
\]
\end{lemma}
Lov\'asz and Schrijver (1991) define an operator $N_0$ by:
\[ N_0(K)= \bigcap_{i=1,\cdots,d} \left\{\left(K \cap \{x:  x_i = 0
\}\right) + \left(K \cap \{x: x_i = x_0 \}
\right)\right\}.\] Thus, $N(K)\subseteq N_0(K)$.
The iterated operator $N_0^r$, $N_0$-rank of inequalities,
polytopes and convex cones are defined analogously to the
corresponding definitions of $N$-and $N_+$-ranks.

Lemma \ref{lem:2.1} shows that an inequality will be valid for $N(K)$
if it is valid for $K \cap \{x: x_i = 0 \}$ and $K \cap \{x: x_i = x_0
\}$ for some $i$. In order to iterate Lemma \ref{lem:2.1}, we first
need the following lemma. It is stated in terms of the faces of $Q$,
which can be obtained by intersecting $Q$ with hyperplanes of the form
$\{x: x_i=0\}$ or $\{x: x_i=x_0\}$. Similar insights for a procedure
related to the $N$- procedure were discussed by Balas (1974).

\begin{lemma}
\label{lem:face} Let $F$ be any face of $Q$. Then 
\[ N\left(K \cap F\right)
=  N(K) \cap F.
\] Similarly for $N_+$ and $N_0$. 
\end{lemma}

\proof 
``$\subseteq$'' is clear from the definitions.  For the
converse, let $x\in N(K)\cap F$. This means that there exists a matrix
$Y\in M(K)$ with $Ye_0=x$. Since $Ye_i\in K\subseteq Q$ and
$Y(e_0-e_i)\in K\subseteq Q$ and their sum $Ye_i + Y(e_0-e_i)=Ye_0$
belongs to the face $F$ of $Q$, we have that $Ye_i$ and $Y(e_0-e_i)$
must belong to $F$, by definition of a face. Thus, $Ye_i\in K\cap F$
and $Y(e_0-e_i)\in K\cap F$ for all $i$ implying that $Y\in M(K\cap
F)$ and $x\in N(K\cap F)$. The proof for $N_+$ is identical.
\qed

Iterating Lemma \ref{lem:face}, we get:
\begin{corollary}
\label{cor:face} Let $F$ be any face of $Q$. Then, for any $r$, 
\[ N^r\left(K \cap F\right)
=  N^r(K) \cap F.\] 
Similarly for $N_+$ and $N_0$. 
\end{corollary}

Repeatedly using Lemma \ref{lem:2.1} and Lemma \ref{lem:face} (or
Corollary \ref{cor:face}), we can derive a condition that an
inequality be valid for $N^r(K)$. This, in particular, proves Theorem
\ref{thm:LS}.

\begin{theorem} \label{rankn}
 $N^r_+(K)\subseteq N^r(K) \subseteq  N^r_0(K) \subseteq
\tilde{N}^r_{0}(K)$ where 
\[ \tilde{N}^r_0(K) = \bigcap_{\{J\subseteq \{1,\cdots,d\}:
|J|=r\}} \sum_{\{(J_0,J_1) \mbox{ partitions
of } J\}} \left(K \cap \{x:
x_i=0 \mbox{ for }i\in J_0 \mbox{ and } x_i=x_0 \mbox{ for }  
i\in J_1\}\right).
\]
\end{theorem}
We should point out that even though
$N_0(K) = \tilde{N}_0(K)$ and
$N^d_+(K) = N^d(K) = N^d_0(K) = \tilde{N}^d_{0}(K)$,
$N_0^r(K)$ is not necessarily equal to
$\tilde{N}^r_{0}(K)$, if $2 \leq r \leq (d-1)$.
For example, for $K=\{x\in Q: x_1+x_2+x_3\leq 1.5
x_0\}$, one can show that $(1,0.5,0.5,0.5)\in
(\tilde{N}_0^2(K)\setminus N_0^2(K))$.

For $N_+(K)$, Lov\'asz and Schrijver (1991) give a different
condition for the validity of an inequality.  In the statement
of the next lemma, the
assumption that $a\geq 0$ is without loss of generality (by flipping
coordinates if necessary, as shown in Proposition \ref{prop:2.1}).

\begin{lemma}[Lov\'asz and Schrijver (1991)]
\label{lem:2.2}
Let $a \geq 0$. Then
$a^T x \leq \alpha x_0$ is valid for\\
$\left(K \cap \{x: x_i = x_0 \}\right)$
for all
$i$ such that $a_i > 0$, implies
$a^T x \leq \alpha x_0$ is valid for $N_+(K)$.
\end{lemma}

As mentioned previously, the result that clique, odd hole, odd
antihole, odd wheel inequalities for the stable set problem have
$N_+$-rank 1 follows from the above lemma. For the stable set problem
(as for many combinatorial optimization problems), there exists
several important constructions to derive facet-defining valid
inequalities from other facet-defining inequalities. The simplest is
{\it cloning} a clique at a vertex $v$, which consists of replacing
the vertex by a clique, replacing all the edges incident to $v$ by
corresponding edges incident to all clique vertices and substituting
in the inequality the variable for $v$ by the sum of the variables of
the clique vertices. It can easily be shown that the resulting
inequality is valid and facet-defining if the original inequality was
a non-trivial (i.e.\ different from the nonnegativity constraints)
facet-defining inequality. In general, it is not clear how cloning
influences the $N_+$-rank of an inequality. However, if we perform
cloning at the center vertex of an odd wheel inequality, Lemma
\ref{lem:2.2} implies that the $N_+$-rank still remains equal to 1.  
If we perform cloning at one or several vertices of an odd wheel, odd
hole or odd antihole inequality, Lemma \ref{lem:2.2} implies that the
$N_+$-rank is at most 2. Indeed, if we fix any variable (of the
corresponding subgraph) to 1, the resulting inequality can be seen to
be a linear combination of clique inequalities and hence valid for
$N_+(FRAC)$. 

Lemma \ref{lem:2.2} can be extended to derive conditions under which the
$N_+$-rank of an inequality is at most $r$. 

\begin{theorem}
\label{thm:3.5}
Let $a \geq 0$ and let $I_+=\{i: a_i>0\}$.  If $a^T x \leq \alpha x_0$
is valid for\\
$\left(K \cap \{x: x_i = x_0,  \mbox{ for all }
i \in I \}\right)$ for all sets $I\subseteq I_+$ satisfying either of
the following two conditions
\begin{enumerate} 
\item \label{it:1} $|I| = r$,
\item \label{it:2}
 $|I| \leq (r-1)$ and $\sum_{i \in I} a_i > \alpha$,
\end{enumerate}
 then $a^T x \leq \alpha x_0$ is valid for $N^r_+(K)$.
\end{theorem}

Observe, however, that
the result mentioned previously regarding cloning does not follow from
Theorem \ref{thm:3.5}. 

\proof We proceed by induction on $r$. For $r=1$, the result is Lemma
\ref{lem:2.2}. 

Assume now that $r>1$, that the theorem was proved for $(r-1)$ (and for any
inequality and for any convex set $K$), and that the hypothesis is
satisfied for the inequality $a^Tx\leq \alpha x_0$ and $r$. {From}
Corollary \ref{cor:face} and Lemma \ref{lem:2.2}, we know that
$a^Tx\leq \alpha x_0$ is valid for $N_+^r(K)=N_+(N_+^{r-1}(K))$ if it
is valid for $N_+^{r-1}(K) \cap \{x:
x_i=x_0\}=N_+^{r-1}\left(K\cap\{x: x_i=x_0\}\right) $ for all $i\in
I_+$.  This is equivalent to showing that $a^Tx -a_ix_i \leq
(\alpha-a_i) x_0$ is valid for $N_+^{r-1}\left(K\cap\{x:
x_i=x_0\}\right)$.

Now there are two
cases. If $\alpha-a_i<0$ then condition \ref{it:2} implies that $K\cap
\{x: x_i=x_0\}=\emptyset$ and thus any inequality is valid for
$N_+^{r-1}\left(K\cap\{x: x_i=x_0\}\right)=\emptyset$. On the other hand, if
$\alpha-a_i\geq 0$, we can use induction to prove the result. Indeed,
conditions \ref{it:1} and \ref{it:2} for inequality $a^Tx\leq \alpha
x_0$ and $r$ imply that conditions  \ref{it:1} and \ref{it:2} are
satisfied for the inequality $a^Tx-a_i x_i\leq (\alpha-a_i)
x_0$ for $r-1$. Thus, by the inductive hypothesis,  $a^Tx -a_ix_i \leq
(\alpha-a_i) x_0$ is valid for $N_+^{r-1}\left(K\cap\{x:
x_i=x_0\}\right)$, proving the inductive statement.
\qed

For the stable set problem, the above theorem implies that the
$N_+$-rank of a graph is at most its stability number $\alpha(G)$, the
cardinality of the largest stable set in $G$; this was proved in Corollary
2.19 of Lov\'asz and Schrijver (1991). More generally, if we consider
a polytope $P$ for which $P_I$ is only described by inequalities of
the form $a^Tx\leq \alpha x_0$ with $a\geq 0$ (i.e.~it is lower
comprehensive, see Section \ref{sec:lowercomp}) then its $N_+$-rank is
upper bounded by the maximum number of variables that can be set to 1
in $P$ to obtain a unique integral point of $P_I$ (in which the other
variables are thus set to 0).  Similar, more complex, statements can
be made if the polytope is not lower comprehensive.

\subsection{Example 1: Matching polytope}

Consider the complete undirected graph on the vertex set $V$; let $E$
denote its edge set. Let 
$$P := \{x \in {\mathbb R}^E: x(\delta(v)) \leq 1, \forall v \in V,
\,\, 0\leq x \leq e\}.$$ 
In the above, $\delta(v)$ is the set of edges in $E$ that are incident
on $v$; for $S \subseteq E$, $x(S)$ represents $\sum\limits_{j \in S}
x_j$. For $S \subseteq V$, let $E(S)$ refer to the set of edges 
with both endpoints in $S$.  Then the {\em matching polytope}
for the complete graph is
$$P_I := \mbox{conv}\left\{P \cap \{0,1\}^E\right\}.$$ Edmonds
(1965) proved that
$$P_I = \left\{ x \in P:
x(E(S)) \leq \frac{|S|-1}{2} \mbox{ for all $S \subseteq V$
such that $|S|$ is odd}\right\}.$$ The above inequalities are known as
the blossom inequalities. 

\begin{theorem}[Stephen and Tun\c{c}el (1999)]
\label{thm:3.1}
The 
$N_+$-rank of the inequality
\[
x\left(E(S)\right) \leq \frac{|S|-1}{2}
\]
with respect to $P$ is $\frac{|S|-1}{2}$.
\end{theorem}

The fact that the $N_+$-rank is at most $\frac{|S|-1}{2}$ also follows
directly from Theorem \ref{thm:3.5}. Observe that since $d$ is $|V|
(|V|-1)/2$, we derive that the $N_+$-rank of
$P$ is equal to $(\sqrt{1+8d}-1)/4$ if $|V|$ odd and,
$(\sqrt{1+8d}-3)/4$ if $|V|$ even.

{From} Theorem \ref{thm:3.1}, the $N$-rank  of the blossom inequality on
$S$ is at least $\frac{|S|-1}{2}$. Furthermore, using Theorem
\ref{rankn} with $J$ being the complement of a complete bipartite
graph on $\frac{|S|-1}{2}$ and $\frac{|S|+1}{2}$ vertices on each
side, we derive that the $N_0$-rank of a blossom inequality is equal to
$\frac{(|S|-1)^2}{4}$. This uses the fact that $P$ is an integral
polytope if and only if the underlying graph is bipartite. Thus, the
$N$-rank is at most $\frac{(|S|-1)^2}{4}$. These
bounds are to be compared with those derived from Corollary 2.8 of
Lov\'asz and Schrijver (1991) (since a matching in a graph can be viewed as
a stable set in its line graph). Their results imply a lower bound of
$(|S|-2)$ and an upper bound of $\frac{1}{2} (|S|-1)^2-1$. 

\subsection{Example 2} \label{ex2}

Consider
\[
K:= \left\{\pmatrix{x_0 \cr x} \in {\mathbb R}^{d+1}: 
x(S) \,\, \leq \,\, \frac{d}{2} x_0, 
\mbox{ for all } S \subset \{1,2, \ldots, d\}
\mbox{ such that } |S|= \frac{d}{2} + 1,
\,\, 0 \leq x \leq x_0 e
\right\}.
\]
Then
\[
K_I = 
\left\{\pmatrix{x_0 \cr x} \in {\mathbb R}^{d+1}:
\,\, \sum_{i=1}^d x_i \,\, \leq \,\,
\frac{d}{2} x_0,
\,\, 0 \leq x \leq x_0 e\right\}.
\]

Theorem \ref{rankn} implies that the $N$-rank of $\sum_{i=1}^d x_i
\leq \frac{d}{2}x_0$ is at most $(d-2)$, while Theorem \ref{thm:3.5}
implies that the $N_+$-rank is at most $\frac{d}{2}$. 
These bounds are actually attained and this is discussed in Section
\ref{details}. We also show in that section that 
the positive semidefiniteness constraint does not help for many
iterations.

\section{Lower bounds on the $N$- or $N_+$-rank} \label{sec:low}

In this section, we provide lower bounds on the $N$- and
$N_+$-rank. We also show a situation in which the positive semidefiniteness
constraints do not help at all and both the $N$-rank and the
$N_+$-rank of a polytope is $d$.

We first provide a way to derive points in $N_+(P)$ in certain cases. 
For $x \in {\mathbb R}^d$ define
\[
x_i^{(j)}  :=\left\{\begin{array}{ll}
x_i & \mbox{ if } i \neq j;\\
0 & \mbox{ if } i =j.
\end{array}
\right.
\]
So, $x^{(j)} = x - x_j e_j$. Throughout this section, let
$K=\left\{\pmatrix{\lambda \cr \lambda x}: x\in P, \lambda \geq 0\right\}$.

\begin{theorem}
\label{thm:3.2}
Let $\bar{x} \in P$ such that
\[
\bar{x}^{(j)} \mbox{ and }
(\bar{x}^{(j)} + e_j) \in  P,
\mbox{ for all } j 
\mbox{ such that } 0 < \bar{x}_j < 1.
\]
Then $x\in N_+(P)$.
\end{theorem}
Simply stated, this result says that if we can replace any coordinate
of $x$ (strictly between 0 and 1) by 0 and 1 and remain in $P$
then $x\in N_+(P)$. 

\proof
We define
\[
Y(x)  := 
\pmatrix{1 \cr x} \pmatrix{1, x^T} +
\mbox{Diag}\pmatrix{0 \cr
x_1 - x_1^2 \cr
x_2 - x_2^2 \cr
\vdots \cr
x_d -x_d^2}.
\]
By definition, $Y(\bar{x}) \in {\cal S}^{d+1}$,
$Y(\bar{x})e_0 = \mbox{diag}\left(Y(\bar{x})\right)
= \pmatrix{1 \cr \bar{x}} \in K$.
Moreover,
\[
Y(\bar{x}) e_j = \bar{x}_j\pmatrix{1 \cr \bar{x}^{(j)} + e_j},
\mbox{ for all } j \in \{1,2, \ldots, d\};
\]
therefore, $Y(\bar{x}) e_j \in K$ for all
$j \in \{1,2, \ldots, d\}$. Similarly,
\[
Y(\bar{x}) (e_0 - e_j) = (1- \bar{x}_j)\pmatrix{ 1 \cr \bar{x}^{(j)}},
\mbox{ for all } j \in \{1,2, \ldots, d\};
\]
therefore, $Y(\bar{x})(e_0 - e_j) \in K$ for all
$j \in \{1,2, \ldots, d\}$.
Finally, since
\[
\mbox{Diag}\pmatrix{0 \cr
x_1 - x_1^2 \cr
x_2 - x_2^2 \cr
\vdots \cr
x_d -x_d^2} \succeq 0
\mbox{ and }
\pmatrix{1 \cr x} \pmatrix{1, x^T}
\succeq 0,
\]
for all $0 \leq x \leq e$,
we have $Y(\bar{x}) \succeq 0.$
Therefore,
$Y(\bar{x}) \in M_+(P)$ and
$x \in N_+(P)$
as desired.
\qed

\begin{figure}[htbp]
\hspace*{\fill} \epsfig{figure=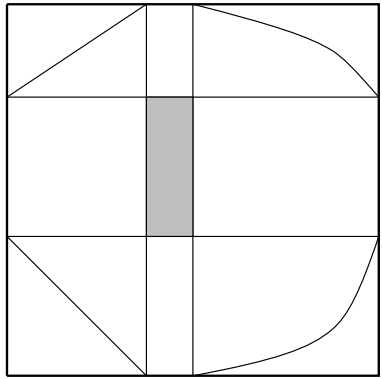,width=2in} \hspace*{\fill}

\caption{\label{fig2} Convex set satisfying the condition of Corollary
\protect\ref{corequal}.}
\end{figure}

As a corollary, we derive the following (see Figure \ref{fig2} for an
illustration).
\begin{corollary} \label{corequal}
Let $P$ be such that $(P\cap\{x: x_j=0\})+ e_j = P\cap
\{x:x_j=1\}$ for all $j\in\{1,\cdots,d\}$. Then $$N_+(P)=N(P)=N_0(P)=
\bigcap_{j\in\{1,\cdots,d\}} \{x: x^{(j)}\in P\}.$$
\end{corollary}

\proof
Let $C=\bigcap_{j\in\{1,\cdots,d\}} \{x: x^{(j)}\in P\}$. 
By Lemma \ref{lem:2.1}, we know that $N_+(P)\subseteq N(P)\subseteq
N_0(P) \subseteq 
C$. On the other hand, Theorem \ref{thm:3.2} shows that $C\subseteq
N_+(P)$. 
\qed

In the proof of Theorem \ref{thm:3.2}, we constructed
a $Y \in M_+(P)$ such that a certain $x \in P$ would also be in $N_+(P)$.
The idea of the proof suggests a stronger technique to achieve such a goal.
We define
\[
Y(x) := \pmatrix{1 \cr x} \pmatrix{1, & x^T}
+ \mbox{Diag}\pmatrix{0 \cr x_1 - x_1^2 \cr
\vdots \cr x_d -x_d^2}
+ \pmatrix{ 0 & 0^T \cr 0 & B(x)},
\]
where $B(x) \in {\mathcal S}^d$, $\mbox{diag}(B) = 0.$
Then clearly we have $Y(x) \in {\mathcal S}^{d+1},$ $Y(x) e_0 = \mbox{diag}
\left(Y(x)\right).$
Moreover, using the Schur complement of $(Y(x))_{00}$ in $Y(x)$, we have
\[
Y(x) \succeq 0 \mbox{ iff } B(x) + \mbox{Diag}\pmatrix{x_1 - x_1^2 \cr
\vdots \cr x_d -x_d^2} \succeq 0.
\]
The latter can be assured in many simple ways, for example by
diagonal dominance:
It suffices to choose $B_{ij}$ such that
\[
|B_{ij}| \leq \frac{1}{2}
\min\left\{ \frac{x_i-x_i^2}{\mbox{\# of nonzeros in column $i$}},
\frac{x_j-x_j^2}{\mbox{\# of nonzeros in column $j$}} \right\}.
\]

The entries of such a $B(x)$ will be further restricted
by the condition
$Y(x)e_i \in K$ for every $i \in \{1,2, \ldots, d\}$ and
$Y(x)(e_0-e_i) \in K$ for every
$i \in \{1,2, \ldots,d\}$.
If this condition is verified for some $B(x)$ then the
above argument would imply
$x \in N_+(P)$. In the case of Theorem \ref{thm:3.2}, we
utilized diagonal dominance; because of the
special structure of $P$,
we could choose $B(x) := 0$ and satisfy all the conditions for
$x \in N_+(P)$.

\subsection{Example 3: Infeasibility detection}

We now give an example where both $N$ and $N_+$ require $d$
iterations, showing that Theorem \ref{thm:LS} cannot be improved. This
result was independently obtained by Cook and Dash (1999) who also
show additional results regarding the rank of
inequalities. Previously, the worse example known in terms the number
of repeated $N_+$ iterations needed to obtain $K_I$ was the matching
polytope results of Stephen and Tun\c{c}el (1999) where the
$N_+$-rank was of the order of $\sqrt{d}$. 

Let \[
P (p):= \left\{ x \in {\mathbb R}^d: \left\|x-\frac{1}{2} e
\right\|_1 \leq \frac{p}{2}
\right\}.
\]

\begin{theorem}
\label{thm:3.6}
For $0<p < d$, $N_+(P(p)) \supseteq P(p-1)$. Furthermore, $P(1)\neq
\emptyset$ while $P_I(d-1)=\emptyset$.  Thus, the $N_+$ procedure requires
$d$ iterations to prove $P_I(d-1) =\emptyset.$
\end{theorem}

\proof
Follows from Corollary \ref{corequal}. (In fact this corollary
characterizes precisely $N_+(P(p))$.)
\qed

One interesting feature of the example above is that $P(d-1)$ can be
described by $2^d$ inequalities, contains no integral point, but no
inequality can be removed without creating an integral point. This is
actually an extreme situation in this regard as shown by the following
result of Doignon (1973). 
Suppose we are given a set of $m$ linear inequalities
\[
a_i^T x \leq  b_i, 
\mbox{ for all } i \in J,
\]
where $x \in {\mathbb R}^d$ and $|J| \geq 2^d$.  A theorem of Doignon
(1973) implies that if this system does not contain any
integer points then there is a subsystem (of this system) with at most
$2^d$ inequalities which does not have an integer solution. 
Doignon's Theorem is an integer analog of Helly's Theorem.  

\subsection{Example 2, continued} \label{details}

In Section \ref{ex2}, we have shown that the $N$-rank and the
$N_+$-rank of 
\[
K:= \left\{\pmatrix{x_0 \cr x} \in {\mathbb R}^{d+1}: 
x(S) \,\, \leq \,\, \frac{d}{2} x_0, 
\mbox{ for all } S \subset \{1,2, \ldots, d\}
\mbox{ such that } |S|= \frac{d}{2} + 1,
\,\, 0 \leq x \leq x_0 e
\right\},
\]
are at most $(d-2)$ and $d/2$, respectively. Here we claim that these
bounds are attained. 

\begin{theorem}
\label{thm:3.8}
The $N$-rank of $\sum_{i=1}^d x_i \leq \frac{d}{2}$
relative to $K$ is $(d-2)$.
The $N_+$-rank of the same inequality relative to
$K$ is $\frac{d}{2}.$

oreover, for $r \leq \frac{d}{2} - \sqrt{d}+3/2$,
the optimum values of
\[
\max\{e^Tx: \,\, x \in N^r(K) \}
\mbox{ and }
\max\{e^Tx: \,\, x \in N^r_+(K) \}
\]
are the same. 
\end{theorem}

Our proof of the first statement of the theorem, saying that the
$N$-rank is $(d-2)$ is lengthy and is not included here. The proof of
the remainder of the theorem appears partly in this section and partly
in the Appendix.
The theorem indicates that the positive semidefiniteness constraint does not
help for $(d/2 - o(d))$ iterations.

Unfortunately, neither Theorem \ref{thm:3.2} nor Corollary
\ref{corequal} is useful here. Instead, exploiting the symmetry (and
convexity of $N(K)$ and $N_+(K)$), we will only consider points in
$N^r(K)$ or $N_+^r(K)$ such that $x_i$ takes only three possible
values, $0$, $1$ and a constant $\alpha$. Letting $n_0$ denote the
number of $x_i$ set to 0 and letting $n_1$ denote the number of $x_i$
set to 1, we define $c(r,n_0,n_1)$ to be the largest common value
$\alpha$ of the remaining $(d-n_0-n_1)$ coordinates of $x$ such that
$x\in N^r(K)$. We define $c_+(r,n_0,n_1)$ similarly with respect to
$N_+^r(K)$.

By symmetry, such a point $x$ belongs to $N^r(K)$ (resp.\ to $N^r_+(K)$) if
there exists a symmetric matrix $Y\in M(K)$ (resp.\ $Y\in M_+(K)$) of the form
\[
Y(n_0, n_1;\alpha, \beta) := \pmatrix{1 & e^T & 0 & \alpha e^T \cr
e & e e^T & 0 & \alpha e e^T \cr
0 & 0 & 0 & 0 \cr
\alpha e & \alpha e e^T & 0 & (\alpha -\beta)I + \beta e e^T},
\]
for some value $\beta$; here the columns of $Y$ are partitioned in
the way that the first column corresponds to the homogenizing
variable $x_0$, the next $n_1$ columns correspond to those $x_j$ that
are set to one, the next $n_0$ columns correspond to those $x_j$ set
to zero and the remaining $(d - n_0- n_1)$ columns correspond to the
remaining $x_j$'s (which are set to $\alpha$).

For $r=0$ and $n_1\leq d/2$, we see by plugging $x$ into the
description of $K$ that
\begin{equation}\label{basecase}
c(0,n_0,n_1)=c_+(0,n_0,n_1)=\left\{ \begin{array}{ll}
\frac{d/2-n_1}{d/2+1-n_1} & \mbox{if } n_0\leq d/2-1, \\ 1 &
\mbox{otherwise.} \end{array} \right.
\end{equation}

For $r>0$, the condition that $Y\in M^r(K)$ is equivalent to
$\frac{\beta}{\alpha} \leq c(r-1,n_0,n_1+1)$ (corresponding to
$Ye_i\in M^{r-1}(K)$) and $\frac{\alpha-\beta}{1-\alpha} \leq
c(r-1,n_0+1,n_1)$ (corresponding to $Y(e_0-e_i)\in
M^{r-1}(K)$). Eliminating $\beta$, we derive: $$ c(r,n_0,n_1)=
\frac{c(r-1,n_0+1,n_1)}{1-c(r-1,n_0,n_1+1)+c(r-1,n_0+1,n_1)}.$$
The condition that $Y\succeq 0$ reduces to (by taking a Schur complement) 
$(\alpha-\beta)I+(\beta-\alpha^2)ee^T\succeq 0$ (where the matrices
have size $(d-n_0-n_1)\times (d-n_0-n_1)$, or $\alpha-\beta\geq
0$ and $\alpha-\beta+(d-n_0-n_1)(\beta-\alpha^2)\geq 0$. This can be
seen to imply 
that 
\begin{eqnarray*} c_+(r,n_0,n_1)& = & \min\left(
\frac{c_+(r-1,n_0+1,n_1)}{1-c_+(r-1,n_0,n_1+1)+c_+(r-1,n_0+1,n_1)}, 
\right. \\ &
& \left.
\frac{(d-n_0-n_1-1) c_+(r-1,n_0,n_1+1)+1}{d-n_0-n_1}\right).
\end{eqnarray*}

Observe that the $N$-rank (resp.\ the $N_+$-rank) of $K$ is the
smallest integer $r$ such that $c(r,0,0)=\frac{1}{2}$ (resp.\
$c_+(r,0,0)=\frac{1}{2}$). 
Theorem \ref{thm:3.8} hence follows from the following proposition. 

\begin{proposition} \label{propapp}
~
\begin{enumerate}
\item \label{firstitem}
$ c(d-3,0,0)=\left\{\begin{array}{ll} \frac{1}{2}+\frac{1}{5d-6} &
\mbox{ if $d$ is even} \\ \frac{1}{2}+\frac{1}{10d-20} & \mbox{ if $d$
is odd} \end{array} \right.$,
\item
$c_+(d/2-1,0,0)>0.5$,
\item
For any $r, n_0,n_1$ such that $r+n_0+n_1\leq d/2-\sqrt{d}+3/2$, we
have $c(r,n_0,n_1)=c_+(r,n_0,n_1)$. 
\end{enumerate}
\end{proposition}

The proof of \ref{firstitem} is obtained by solving explicitly the
recurrence for $c$; the details however, are omitted. The proof of the
rest of the proposition is given in the Appendix. 

\begin{figure}[htbp]
\hspace*{\fill}  \epsfig{figure=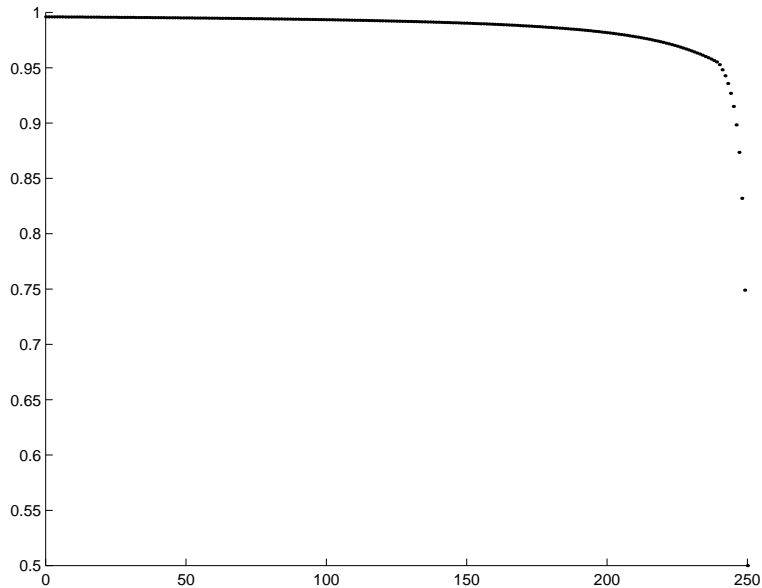,width=4in} \hspace*{\fill}

\caption{\label{fig:nplus} Plot of $c_+(r,0,0)$ for $d=500$ as a
function of $r$.}
\end{figure}

Theorem \ref{lower} in the Appendix actually illustrates a peculiar
behavior of the $N_+$ operator (as well as the $N$ operator) on
this example. In
cutting plane procedures, it is usual that the improvement due to the
addition of a cutting plane (or a batch of them) decreases as the algorithm
progresses. However, Theorem \ref{lower} shows that $$\max\{e^Tx: x\in
N^r_+(K)\} = d c_+(r,0,0)>d \left(1 -
\frac{1}{d/2+1-r}\right).$$ Hence, as illustrated on Figure
\ref{fig:nplus} for $d=500$, the improvement in objective function
value is negligible for many iterations and only towards the end
increases considerably. We should point out, however, that the
procedures $N$ and $N_+$ are such that the number of ``important''
inequalities generated in each iteration could potentially increase
tremendously in later iterations.

\section{Additional properties} \label{sec:lowercomp}

A nonempty convex set $P \subseteq {\mathbb R}^d_+$ is called
{\em lower comprehensive} if for every $x \in P$,
every $y \in {\mathbb R}^d_+$ such that $y \leq x$
is also in $P$.

\begin{definition}
\label{defi:3.1}
Let $v \in \{0,1\}^d$.  A convex set $P \subseteq [0,1]^d$ is said to
be {\em a convex corner with respect to $v$} if there exists a linear
transformation $A$ of $\{0,1\}^d$ onto itself such that $Av =0$ and
$AP$ is lower comprehensive.
\end{definition}

\begin{theorem}
\label{thm:3.3}
If $P$ is a convex corner  with respect to
$v \in \{0,1\}^d$ then so are
$N(P)$ and $N_+(P)$.
\end{theorem}

\proof
By Proposition \ref{prop:2.1} and the definitions,
it suffices to prove that if $P$ is lower comprehensive
then so are $N(P)$ and $N_+(P)$. Let $P$
be lower comprehensive and $x \in N(P)$.
It suffices to show that $(x - x_j e_j) \in N(P)$ for
every $j$ such that $x_j > 0$.
Without loss of generality suppose
$j = 1$ and $x_j > 0$.
Then there exists $Y \in M(P)$ such that $Y e_0 =
\pmatrix{1 \cr x}$.
Let
\[
\bar{Y}_{ij} \,\, := \,\, \left\{
\begin{array}{cl}
Y_{ij} & \mbox{if $i \neq 1$ or $j \neq 1$};\\
0 & \mbox{otherwise}.
\end{array}
\right.
\]
Then using the fact that $P$ is lower comprehensive, it is easy
to see that $\bar{Y} \in M(P)$.
Since the above argument applies to every $j$
such that $x_j > 0$, we proved that
$N(P)$ is lower comprehensive.
\par
We can prove that $N_+(P)$ is lower comprehensive by a
very similar argument. We only have to note that if
$Y \in M_+(P)$ then the corresponding $\bar{Y}$ constructed
as above will be positive semidefinite (in addition to
satisfying $\bar{Y}e_j \in K$ for every $j \in \{0,1,2,
\ldots, d\}$ and $\bar{Y}(e_0-e_j) \in K$
for every $j \in \{1,2, \ldots,d\}$)
since every principal minor of $\bar{Y}$ is a principal minor
of $Y$ and $Y$ is positive semidefinite.
\qed

A similar fact, in a less general form, was observed independently
by Cook and Dash (1999).


\section{General conditions on the strength of the semidefinite
constraint} \label{sec:gen}

In this section, we derive general conditions under which the
positive semidefiniteness constraint is not useful. This can be expressed in
several ways as
\begin{itemize}
\item
$M(K)=M_+(K)$, or as
\item
$N(K)=N_+(K)$ or even as
\item
$\max\{c^Tx:
x\in N(K)\}=\max\{c^Tx: x\in N_+(K)\}$ for some given $c$. 
\end{itemize}

First, we rewrite condition (ii) of Definition \ref{defi:M}. 
Since $Y$ is symmetric,
\[
u^T Y v \, \geq \, 0, \,\,\, \forall \, u \in Q^*, \, v \in K^*
\iff
u^T Y v + v^T Y u \, \geq \, 0, \,\,\, \forall \, u \in Q^*, \, v \in K^*.
\]
Using the fact that $u^T Y v + v^T Y u  \, = \, \mbox{Tr} \left(
Y(u v^T + v u^T) \right)$, we see that condition $(ii)$ is
also equivalent to 
\beann
(ii)^{''} & &  Y \in \left[T(K)\right]^*,
\eeann
where
\[
T(K) \,\, := \,\, \mbox{cone} \left\{
u v^T + v u^T: u \in Q^*, v \in K^* \right\}
\,\, = \,\, \mbox{cone} \left\{
u v^T + v u^T: u \in \mbox{ext}(Q^*), v \in \mbox{ext}(K^*)
\right\}.
\]
Let's define
\[
D := \left\{ Y \in {\cal S}^{d+1} : \mbox{diag} (Y) = Y e_0
\right\}.
\]
Note that the cone (more specifically,
the subspace in this case) dual to $D$ in the space ${\cal S}^{d+1}$
is the orthogonal complement
of $D$.
\[
D^*  = D^{\perp} = \left\{ \sum\limits_{i=1}^{d}
\alpha_i ( E_{ii} - E_{0i}) : \alpha \in {\mathbb R}^d \right\},
\]
where $E_{ij} := e_ie_j^T + e_j e_i^T$.
We have

\begin{theorem}
\label{thm:M=M+}
\[
M_+(K) = M(K) \,\,\,\, \mbox{ if and only if } \,\,\,\,
T(K) + D^{\perp} \supseteq {\cal S}_+^{d+1}.
\]
\end{theorem}

\proof
By definition of the sets $M(K)$, $M_+(K)$, we have
\[
M(K) = M_+(K) \iff \left[T(K)\right]^* \cap D
= \left[T(K)\right]^* \cap D \cap {\cal S}_+^{d+1}.
\]
Since the inclusion $\left[T(K)\right]^* \cap D \supseteq
\left[T(K)\right]^* \cap D \cap {\cal S}_+^{d+1}$ is clear, we have
\[
M(K) = M_+(K) \iff \left[T(K)\right]^* \cap D
\subseteq {\cal S}_+^{d+1}.
\]
Noting that
\[
\left[T(K)\right]^* \cap D
\subseteq {\cal S}_+^{d+1} \,\, \iff
\,\, (\left[T(K)\right]^* \cap D)^* \supseteq {\cal S}_+^{d+1},
\]
(we used the fact that ${\cal S}^{d+1}_+$
is self dual under the trace inner-product, in the space
${\cal S}^{d+1}$) and that 
\[
(\left[T(K)\right]^* \cap D)^* = T(K) +D^*,
\]
we conclude
\[
M(K) = M_+(K) \,\,\,\, \mbox{ if and only if } \,\,\,\,
T(K)+D^{\perp} \supseteq {\cal S}_+^{d+1}.
\]
\qed

This theorem completely characterizes when $M$ and $M_+$ differ or are
equal. To make the condition more easily tractable, we can give a more
explicit description of $T(K)+D^{\perp}$.  Define $F(K)$ to be set of
all $v = \pmatrix{v_0 \\ \cr {\bar v} \cr} \in {\mathbb R}^{d+1}$ such
that $-{\bar v}^T x \leq v_0$ is a facet of $P$ (or, more generally,
for non-polyhedral convex sets, $F(K)$ describes a set of valid
inequalities exactly characterizing $P$).  Note that $F(K)$ can be
taken as the set of extreme rays of $K^*$.  We arrive at the identity
\beann
T(K)+D^{\perp} & = &
\mbox{cone} \left\{(e_i v^T + v e_i^T), \,\,
i \in \{1,2, \ldots, d\}, v \in F(K);\right.\\
& & \,\,\,\,\,\,\,\,\,\,\,\,\,\,\,\,\,\,
\left[(e_0 -e_i) v^T + v(e_0 -e_i)^T\right], \,\, i \in \{1, 2, \ldots, d\},
v \in F(K);\\
& & \,\,\,\,\,\,\,\,\,\,\,\,\,\,\,\,\,\,
\left. (E_{ii} - E_{0i}), \,\, i \in \{1,2, \ldots, d\},
\right\}
\eeann
where we have used the fact that $E_{0i}-E_{ii}\in T(K)$ since $e_i\in F(K)$. 
So, $M_+(K) = M(K)$ iff for every $x \in {\mathbb R}^{d+1}$,
we can express $xx^T$ as an element of the above cone
$(T(K)+D^{\perp})$.

Consider the clique on four vertices and the corresponding
LP relaxation $FRAC$ of the stable set problem (with the edge and
nonnegativity constraints only). For this example,
\[
Y:=\pmatrix{ 1 & \frac{1}{3} & \frac{1}{3} &\frac{1}{3} & \frac{1}{3} \cr
\frac{1}{3} & \frac{1}{3} & 0 & 0 & 0 \cr
\frac{1}{3} &0 & \frac{1}{3} & 0 & 0 \cr
\frac{1}{3} &0 & 0 & \frac{1}{3} & 0 \cr
\frac{1}{3} &0 & 0 & 0 & \frac{1}{3}} \in 
\left[T(K)\right]^* \cap D;
\]
but clearly $Y \notin {\mathcal S}^{d+1}_+$. A proof of this is provided
by the incidence vector of the clique inequality on the four vertices:
\[
(e_0-e_1-e_2-e_3-e_4)^T Y (e_0-e_1-e_2-e_3-e_4) = - \frac{1}{3}.
\]
This means, for $x:= (e_0-e_1-e_2-e_3-e_4)$,
$xx^T$ is not in the convex cone $\left(T(K) + D^{\perp}\right)$.

Now, we relate these findings to  $N(K)$ and $N_+(K)$.

\begin{corollary}
\label{cor:6.1a}
If $\left(T(K) + D^{\perp}\right)
\supseteq {\mathcal S}^{d+1}_+$,
then
$N_+^r(K) = N^r(K)$ for every $r \geq 0$.
\end{corollary}

\proof
Trivial for $r=0$.
By Theorem \ref{thm:M=M+}, the assumption of the corollary implies
$N_+(K) = N(K)$. By Theorem \ref{thm:LS},
$N(K) \subseteq K$. Thus,
\[
\left(T(N(K)) + D^{\perp} \right)
\supseteq \left(T(K) + D^{\perp}\right) \supseteq {\mathcal S}^{d+1}_+.
\]
Now, applying Theorem
\ref{thm:M=M+} recursively, we obtain the desired result.
\qed

Now we look at the weaker condition that $N(K)=N_+(K)$. 

\begin{theorem}
\label{thm:N=N+}
$N_+(K) = N(K)$ if and only if
for every $s \in {\mathbb R}^{d+1}$,
\[
\mbox{Diag}(s) \in T(K) + D^{\perp} + {\cal S}_+^{d+1}
\,\,\,\,\, \mbox{ implies } \,\,\,\,\,
\mbox{Diag}(s) \in T(K) + D^{\perp}.
\]
\end{theorem}

Before proceeding with the proof, observe that, for any convex cone
${\cal K} \subseteq {\cal S}^{d+1}$, we have
\bea
\label{eqn:adjoint-dual}
\left[\diag({\cal K})\right]^* \,\,
= \,\, \left\{s \in {\mathbb R}^{d+1}: \,\, Diag(s) \in {\cal K}^* \right\}.
\eea

\proof
As in the proof of Theorem \ref{thm:M=M+},
we obtain
\[
N(K) = N_+(K) \iff \mbox{diag}\left(\left[T(K)\right]^* \cap D \right)
\subseteq \mbox{diag}\left(\left[T(K)\right]^* \cap D \cap {\cal S}_+^{d+1}\right).
\]
Using equation (\ref{eqn:adjoint-dual}) and the
proof technique of Theorem \ref{thm:M=M+}, we find
\[
N(K) = N_+(K) \,\,\,\, \mbox{ if and only if }
\]
\[
\mbox{for every } s \in {\mathbb R}^{d+1},
\mbox{Diag}(s) \in T(K) + D^{\perp} + {\cal S}_+^{d+1}
\,\,\,\,\, \mbox{ implies } \,\,\,\,\,
\mbox{Diag}(s) \in T(K) + D^{\perp}.
\]
\qed

We should compare this result to Lemma 1.2 of Lov\'asz and Schrijver (1991).
Note that our result is also based on cone duality, we also
characterize the dual cones of $N(K)$ and $N_+(K)$; but,
we only work in the space of symmetric matrices
instead of the larger space of all matrices. As a result,
the dependence of the characterization
on the skew symmetric matrices is eliminated and
our description is more explicit.

Our ideas in the geometric characterizations above
are also applicable in comparing the weaker procedure
$N_0$ to $N$.
Recall
\[
N_0(K):= \bigcap_{i=1,\cdots,d} \left\{\left(K \cap \{x:  x_i = 0
\}\right) + \left(K \cap \{x: x_i = x_0 \}
\right)\right\}.
\]
We define
\[
M_0(K) := \left\{Y \in {\mathbb R}^{(d+1)\times(d+1)}:
Ye_0 = Y^Te_0 = \mbox{diag}(Y),
u^T Y v \geq 0, \forall u \in Q^*, v \in K^* \right\},
\]
the main difference with $M$ is that $Y$ is not necessarily
symmetric. 
As is mentioned by Lov\'asz and Schrijver (1991),
we have
\[
N_0(K)= \left\{Y e_0: Y \in M_0(K) \right\}.
\]
We further define
\[
T_0(K):=\mbox{cone}\left\{u v^T : u \in Q^*, v \in K^* \right\},
\mbox{ and }
D_0:= \left\{Y \in {\mathbb R}^{(d+1)\times(d+1)}:
Ye_0 = Y^T e_0 = \mbox{diag}(Y) \right\}.
\]
Then
\[
Y \in M_0(K) \mbox{ iff } Y \in \left(\left[T_0(K)\right]^* \cap D_0
\right),
\]
where $\left[T_0(K)\right]^*$ is the dual
of $T_0(K)$ in ${\mathbb R}^{(d+1)\times(d+1)}$
under the trace inner-product.

\begin{theorem}
\label{thm:6.3}
$M_0(K) = M(K)$ iff $\left(T_0(K) + D_0^{\perp}\right)
\supseteq \left\{(e_i e_j^T - e_j e_i^T):
i,j \in \{1,2, \ldots,d\} \right\}.$
\end{theorem}

\proof
As we showed,
$M_0(K) = \left[T_0(K)\right]^* \cap D_0$
and it is clear from the definitions that $M(K) =
\left[T_0(K)\right]^* \cap D_0 \cap {\mathcal S}^{d+1}.$
Note that
\[
D_0^{\perp} = \mbox{span}\left\{
e_ie_i^T-e_0e_i^T, e_i e_i^T - e_i e_0^T:
i \in \{1,2, \ldots, d\}\right\}.
\]
Thus,
\[
\pm (e_0e_i^T-e_ie_0^T) \in \left(T_0(K) + D_0^{\perp}\right),
\forall i \in \{1,2, \ldots,d\}.
\]
Let $\tilde{\mathcal S}^{d+1}$ denote the subspace of $(d+1) \times
(d+1)$ skew-symmetric matrices with real entries.
Therefore,
\[
\left(T_0(K) + D_0^{\perp}\right) \supseteq \tilde{\mathcal S}^{d+1}
\mbox{ iff } \left(T_0(K) + D_0^{\perp}\right) \supseteq
\left\{(e_i e_j^T - e_j e_i^T):i,j \in \{1,2, \ldots,d\} \right\}.
\]
Now, using elementary cone geometry on closed convex cones and the
definitions, we have the following string of equivalences:
\beann
\left(T_0(K) + D_0^{\perp}\right) \supseteq
\left\{(e_i e_j^T - e_j e_i^T):i,j \in \{1,2, \ldots,d\} \right\}
& \mbox{ iff } &
\left(T_0(K) + D_0^{\perp}\right) \supseteq \tilde{\mathcal S}^{d+1}\\
& \mbox{ iff } &
\left[T_0(K)\right]^* \cap D_0 \subseteq {\mathcal S}^{d+1} \\
& \mbox{ iff } &
M_0(K) = M(K).
\eeann
\qed

\begin{corollary}
\label{cor:6.1}
If $\left(T_0(K) + D_0^{\perp}\right)
\supseteq \left\{(e_i e_j^T - e_j e_i^T):
i,j \in \{1,2, \ldots,d\} \right\}$
then
$N_0^r(K) = N^r(K)$ for every $r \geq 0$.
\end{corollary}

\proof
Trivial for $r=0$.
By Theorem \ref{thm:6.3}, the assumption of the corollary implies
$N_0(K) = N(K)$. By Theorem \ref{thm:LS},
$N(K) \subseteq K$. Thus,
\[
T_0(N(K)) \supseteq T_0(K) \supseteq \left\{(e_i e_j^T - e_j e_i^T):
i,j \in \{1,2, \ldots,d\} \right\}.
\]
Now, applying Theorem
\ref{thm:6.3} recursively, we obtain the desired result.
\qed

Let $G$ denote the complete graph on $d$ vertices, and consider the LP
relaxation $FRAC$ of the stable set problem on $G$.
For every $i,j \in \{1,2, \ldots, d\}$ such that $i \neq j$,
we have
\[
(e_0 - e_i -e_j) \in K^* \mbox{ and clearly }
e_i, e_j \in \left(K^* \cap Q^* \right).
\]
Thus, for every $i,j \in \{1,2, \ldots, d\}$ such that $i \neq j$,
we have
\[
e_i (e_0 - e_i -e_j)^T \mbox{ and } e_j e_i^T \in T_0(K),
\mbox{ and } (e_i e_i^T -e_i e_0^T) \in D_0^{\perp}.
\]
This implies,
\[
\left(T_0(K) + D_0^{\perp}\right) \supseteq \left\{(e_i e_j^T - e_j e_i^T):
i,j \in \{1,2, \ldots,d\} \right\}.
\]
Therefore, the condition of Theorem \ref{thm:6.3}
is satisfied and we have $N_0^r(FRAC) = N^r(FRAC)$
for every $r \geq 0$.

As in Theorem \ref{thm:N=N+}, we obtain

\begin{corollary}
\label{cor:6.2}
$N_0(K) = N(K)$ if and only if
for every $s \in {\mathbb R}^{d+1}$,
\[
\mbox{Diag}(s) \in T_0(K) + D_0^{\perp} + \tilde{{\cal S}}^{d+1}
\,\,\,\,\, \mbox{ implies } \,\,\,\,\,
\mbox{Diag}(s) \in T_0(K) + D_0^{\perp}.
\]
\end{corollary}

Instead of comparing $M(K)$ and $M_+(K)$, or $N(K)$ and $N_+(K)$, we
might ask  when are the set of optimal
solutions of both relaxations the same.
This is precisely when
\[
\left[N(K)\right]^* + \pmatrix{-z^* \cr c}
\,\, \supseteq \,\,
\left[N_+(K)\right]^* + \pmatrix{-z^* \cr c},
\]
where $z^*$ is the optimal value of
$\max\left\{c^T x : \,\, \pmatrix{1 \\ \cr x \cr} \in N(K) \right\}$.

Sometimes we are only interested in the bound
provided by the relaxation. This is equivalent to finding
the smallest $z$ for which $\pmatrix{z \\ \cr -c \cr} \in [N(K)]^*$
and the smallest $z^+$ for which
$\pmatrix{z^+ \\ \cr -c \cr} \in [N_+(K)]^*$.

\begin{center}
{\bf References}
\end{center}

\noindent
Balas E. (1974).
Disjunctive programming: Properties of the convex hull
of feasible points. {\em
Management Science Research Report} 348
GSIA, Carnegie Mellon University, Pittsburgh, PA, USA.

\noindent
Balas E., S. Ceria and G. Cornu{\'e}jols (1993).
A lift-and-project
cutting plane algorithm for mixed 0-1 programs.
{\em Math. Prog.} {\bf 58} 295--323.

\noindent
Cook, W., and S.~Dash, private communication, 1999.

\noindent
Doignon, J.-P. (1973). Convexity in cristallographical lattices.
{\em Journal of Geometry} {\bf 3} 71-85.

\noindent
Edmonds, J. (1965).
Maximum matching and a polyhedron with 0,1-vertices.
{\em Journal of Research of the National
Bureau of Standards-B} {\bf 69B} 125--130.

\noindent
Gr\"otschel, M., L.~Lov\'asz and A.~Schrijver (1981). The ellipsoid
method and its consequences in combinatorial optimization. {\em
Combinatorica} {\bf 1(2)} 169--197. 


\noindent
Lov{\'{a}}sz, L. and A. ~Schrijver (1991).  Cones of
matrices and set-functions and 0-1 optimization.  {\em SIAM J. Optim.}
{\bf 1} 166--190.

\noindent
Stephen, T. and L. Tun\c{c}el (1999).  On a
representation of the matching polytope via semidefinite liftings,
{\em Math. Oper. Res.} {\bf 24} 1--7.

\newpage

\begin{appendix}
\section*{APPENDIX}  
\renewcommand{\thesection}{A}
\setcounter{theorem}{0}

In this appendix, we prove Proposition
\ref{propapp} and derive additional properties of $c$ and $c_+$. We
first start with a few preliminary lemmas.

\begin{lemma} \label{lemint}
Assuming $0\leq b<a\leq 1$ and $p>0$, we have
$$a>\frac{a}{1-b+a} >b$$ and
$$a>\min\left(\frac{a}{1-b+a},\frac{(p-1)b+1}{p}\right)>b.$$
\end{lemma}

\proof
First, $a>\frac{a}{1-b+a}$ follows from the fact that $a>0$ and
$a>b$, and this implies also that
$a>\min(\frac{a}{1-b+a},\frac{(p-1)b+1}{p})$. 

On the other hand, we have that $\frac{a}{1-b+a}>b$ iff $a>b-b^2+ab$
iff $(a-b)(1-b)>0$, which follows by assumption.  Furthermore,
$\frac{(p-1)b+1}{p}>b$ iff $(p-1)b+1>pb$ iff $1>b$. As a result, both
terms in the minimum are greater than $b$, and the second part of each
inequality follows.
\qed

This implies the following interlacing property. 

\begin{corollary} \label{cor1}
For any $r\geq 1$ and any $n_0, n_1 \leq \frac{d}{2}-r$, we have
that $$c(r-1,n_0,n_1+1) < c(r,n_0,n_1) < c(r-1,n_0+1,n_1)$$ and 
 $$c_+(r-1,n_0,n_1+1) < c_+(r,n_0,n_1) < c_+(r-1,n_0+1,n_1).$$
\end{corollary}

\proof For $r=1$ and $n_0,n_1\leq \frac{d}{2}-1$, we have that $0\leq
f(r-1,n_0,n_1+1)< f(r-1,n_0+1,n_1) \leq 1$ where $f=c$ or $f=c_+$ by
(\ref{basecase}). Lemma \ref{lemint} now implies the result for $r=1$.

Proceeding by induction on $r$ and assuming true the result for $r-1$,
we derive that $f(r-1,n_0,n_1+1)<f(r-2,n_0+1,n_1+1)<f(r-1,n_0+1,n_1)$,
which implies the result for $r$ by Lemma \ref{lemint}.
\qed

We can now get a lower bound on the coefficients $c$ and
$c_+$. 
\begin{theorem} \label{lower}
For any $r,n_0, n_1$ such that $s=r+n_0+n_1\leq d/2$, we have that 
$$c(r,n_0,n_1)\geq c_+(r,n_0,n_1) > c(0,0,s)=\frac{d/2-s}{d/2+1-s}.$$
In particular, $c_+(d/2-1,0,0)>0.5$.
\end{theorem}

This shows that the $N_+$-rank of $K$ is $d/2$. 

\proof 
For $s\leq d/2$, we have 
$$c(r,n_0,n_1)\geq c_+(r,n_0,n_1) >
c_+(r,0,n_0+n_1)>c_+(0,0,r+n_0+n_1)=c(0,0,s),$$
where we have used Corollary \ref{cor1} twice. 
\qed

\begin{lemma} \label{lemint2}
Let $1\geq  a >b >c\geq  0$ be such that $a-b < b-c$.
Then $$a-b < \frac{a}{1-b+a}-\frac{b}{1-c+b} < b-c.$$
\end{lemma}

\proof
The first inequality is equivalent to 
$$ b\left(\frac{1}{1-c+b}-1\right) <
a\left(\frac{1}{1-b+a}-1\right).$$
This inequality is satisfied since $0< b < a$ and $0 <
\frac{1}{1-c+b}-1 < \frac{1}{1-b+a}-1$ (because
$0< a-b < b-c$). 

For the second inequality, we have that
$$\left(\frac{b}{1-c+b}-c\right)-\left(\frac{a}{1-b+a}-b\right) =
\frac{(b-c)(1-c)}{1-c+b} - \frac{(a-b)(1-b)}{1-b+a}.$$
Moreover, we know that $1-c > 1-b >0$ and
$(b-c)/(1-c+b)> (a-b)/(1-b+a)>0$ since $0< a-b < b-c$.
Multiplying these two inequalities together, we get the desired
inequality. 
\qed

This implies that the coefficients $c(r,n_0,n_1)$ also satisfy the
following differential interlacing property.

\begin{corollary} \label{cor2}
For any $r\geq 1$, any $0\leq n_0\leq d/2-r-2$, any $1\leq n_1\leq
d/2-r$,  we have that
\begin{eqnarray*}
c(r-1,n_0+2,n_1-1)-c(r-1,n_0+1,n_1)& < & c(r,n_0+1,n_1-1)-c(r,n_0,n_1) \\
 &   < & c(r-1,n_0+1,n_1)-c(r-1,n_0,n_1+1).
\end{eqnarray*}
\end{corollary}

\proof For $r=1$, $1\leq n_1\leq \frac{d}{2}-1$ and $n_0\leq d/2-3$, let
$a=c(r-1,n_0+2,n_1-1)$,
$b=c(r-1,n_0+1,n_1)$ and
$c=c(r-1,n_0,n_1+1)$. Observe that
$a=1-\frac{1}{d/2+2-n_1}$, $b=1-\frac{1}{d/2+1-n_1}$ and
$c=1-\frac{1}{d/2-n_1}$, implying that 
$a>b>c$ and
$a-b<b-c$. Thus, Lemma \ref{lemint2} implies the result for $r=1$. 

We now proceed by induction and assume the result true for
$r-1\geq 1$. Defining $a$, $b$ and $c$ as above, we know from Corollary
\ref{cor1} that $a>b>c$ and from the inductive hypothesis that $a-b<
c(r-2,n_0+2,n_1) - c(r-2,n_0+1,n_1+1) < b-c$. Lemma \ref{lemint2} then
implies the result for $r$. 
\qed

Using Corollary \ref{cor2} repeatedly, we derive the following
corollary. 

\begin{corollary} \label{cordiff}
For any $r\geq 1, n_0, n_1\geq 0$ such that $s=r+n_0+n_1\leq d/2$, we
have that
$$c(r-1,n_0+1,n_1)-c(r-1,n_0,n_1+1)<
c(0,1,s-1)-c(0,0,s)=\frac{1}{(d/2+1-s)(d/2+2-s)}.$$ 
\end{corollary}

\proof
Using Corollary \ref{cor2}, we derive
\begin{eqnarray*}
c(r-1,n_0+1,n_1)-c(r-1,n_0,n_1+1) & < & c(r-1,1,n_0+n_1) -
c(r-1,0,n_0+n_1+1) \\ & < & c(0,1,s-1)-c(0,0,s).
\end{eqnarray*}
\qed

\begin{theorem} \label{thmbranch}
For any $r, n_0, n_1\geq 0$ such that $s=r+n_0+n_1\leq d/2-\sqrt{d}+3/2$, we
have that $c(r,n_0,n_1)=c_+(r,n_0,n_1).$
\end{theorem}

\proof
The proof is by induction on $r$. The base case is obvious. Assume the
result is true for $r-1$. This implies that
$c(r-1,n_0+1,n_1)=c_+(r-1,n_0+1,n_1)$ and
$c(r-1,n_0,n_1+1)=c_+(r-1,n_0,n_1+1)$; we denote respectively by $a$
and $b$ these two quantities. The result would then follow if we can
show that
$$\frac{a}{1-b+a}\leq \frac{(p-1) b+1}{p},$$ where $p=d-n_0-n_1$.
This inequality is equivalent to $pa\leq
1-b+a+(p-1)b-(p-1)b^2+(p-1)ab$, or to $(1-b)(a-b)(p-1)\leq 1-b$. Since
$b\leq 1$, we need to prove that $a-b\leq
\frac{1}{p-1}=\frac{1}{d-n_0-n_1-1}.$ 
This follows from Corollary \ref{cordiff} since we have that
$a-b<\frac{1}{(d/2+1-s)(d/2+2-s)}\leq
\frac{1}{(\sqrt{d}-0.5)(\sqrt{d}+0.5)} <\frac{1}{d-1} \leq
\frac{1}{d-n_0-n_1-1}$. 
\qed

\end{appendix}

\end{document}